\documentclass{amsproc}
\usepackage{amsfonts}
\usepackage{amssymb}
\usepackage{amscd}
\usepackage{verbatim}
\usepackage{graphicx}
\usepackage{subfigure}
\usepackage{fancyhdr}
\usepackage[english]{babel}
\usepackage[T1]{fontenc}
\usepackage[latin1]{inputenc}
\usepackage{bm}

\input xy
\xyoption{all}

\newcommand{\mor}[3]{$\xymatrix@1@C=15pt{#3: #1\ar[r]& #2}$}
\newcommand{\appl}[2]{$\xymatrix@1@C=15pt{#1 \ar@{|->}[r]& #2}$}
\newcommand{\pf}{\medskip \noindent {\bf Proof :\ \ }}
\newcommand{\wt}[1]{\widetilde{#1}}
\newcommand{\nwt}[1]{\widetilde{#1}^{\circ}}
\newcommand{\wh}[1]{\widehat{#1}}

\newcommand{\ol}[1]{\overline{#1}}

\newcommand{\cqfd}{\begin{flushright}$\Box$\end{flushright}}
\newcommand{\rad}[1]{{\rm rad}#1}
\newcommand{\BQI}{${\mathcal B}(Q,I)$}

\begin{document}
\bibliographystyle{plain}

\begin{abstract}  We associate to a bound quiver $(Q,I)$ a CW-complex which we denote by \BQI, and call the classifying
space of $(Q,I)$.  We show that the fundamental group of \BQI \ is
isomorphic to the fundamental group of $(Q,I)$. Moreover, we show
that this construction behaves well with respect to coverings. On
the other hand, we study the (co)homology groups of \BQI, and
compare them with the simplicial and the Hochschild (co)homology
groups of the algebra $A=kQ/I$. More precisely, we give sufficient
conditions for these groups to be isomorphic. This generalizes a
theorem due to Gerstenhaber and Schack \cite{GS83}.

\end{abstract}

\title{The classifying space of a bound quiver }
\author[J.~C.~Bustamante ]{Juan Carlos Bustamante
\footnote{
phrases: Fundamental group, bound quivers, simply connected
algebras, schurian algebras, incidence algebras, Hochschild
cohomology, simplicial cohomology of algebras}}
\address{Instituto de Matem\'{a}tica e Estat\'{i}stica,
Universidade de S\H{a}o Paulo, S\H{a}o Paulo, 05508-090, Brazil }
\email{bustaman@ime.usp.br} \maketitle

\thispagestyle{empty}
\section*{Introduction}

Let $A$ be an associative, finite dimensional algebra over an algebraically closed field $k$. It is well-known (see
\cite{BG82}, for instance) that if $A$ is basic and connected, then there exists a connected bound quiver $(Q,I)$ such that
$A \simeq kQ/I$, where $kQ$ is the path algebra of $Q$ and $I$ is an admissible two sided ideal of $kQ$. The pair $(Q,I)$
is then called a {\bf presentation} of $A$. If $Q$ contains no oriented cycles, then $A$ is said to be a {\bf triangular}
$k-$algebra.

For each presentation $(Q,I)$ of $A$, one can define its {\bf fundamental group}, denoted by $\pi_1(Q,I)$ (see \cite{G83, MP83}, for
instance). A triangular $k-$algebra $A$ is said to be {\bf simply connected} if, for every presentation $(Q,I)$ of $A$,  the group
$\pi_1(Q,I)$ is trivial. Simply connected algebras play an important rôle in the representation theory of algebras, since covering
techniques often allow to reduce the study of indecomposable representations of an algebra $A$ to the study of indecomposable
representations of a suitably choosen simply connected algebra \cite{BG82}.

On the other hand, given an algebra $A=kQ/I$, its Hochschild
cohomology groups ${\rm HH}^i(A)$  also give important information
about the simple connectedness, as well as about the rigidity
properties of $A$, see \cite{G64,C91, Skow92}. If moreover $A$
admits a semi-normed basis, one can define the {\bf simplicial
homology}  and { \bf cohomology groups of $A$} with coefficients
in some abelian group $G$, denoted by ${\rm SH}_i(A)$, and ${\rm
SH}^i(A,G)$, respectively see \cite{BG83} (also \cite{MP84,
MdlP99}).

Results in  \cite{AdlP96,BG83, Bus02,PS01, GS83,IZ90} exhibit
several relations  between the groups mentioned above. To a
schurian triangular algebra $A$ is associated a simplicial complex
$|A|$, see \cite{BG83} (and also \cite{Bo84, MP84}). It turns that
in this case the groups $\pi_1(Q,I)$ and the simplicial
(co)homology groups of $A$ are respectively isomorphic to the
fundamental group and the (co)homology groups of $|A|$ see
\cite{Bus02, MdlP99, BG83}.

The main aim of this work is to build a topological space which
would be a geome\-trical model for studying the fundamental groups
and coverings of bound quivers, as well as the simplicial and
Hochschild (co)homology of not necessarily schurian algebras. The
paper is organized as follows.

\medskip
In section $1$, we fix notation and terminology, and recall the definition of the (natural) homotopy relation induced by an
ideal $I$ on the set of walks in a quiver $Q$. We also recall the definition of the fundamental group $\pi_1(Q,I)$.

\medskip We begin section $2$ with a motivating discussion about classifying spaces of small categories, in the sense of $K-$theory
\cite{Se68, Q73}. A particular case of that construction is the simplicial complex $|A|$ associated to an incidence
algebra $A(\Sigma)$. The main part of this section is devoted to the construction of a CW-complex ${\mathcal B}={\mathcal B}(Q,I)$, which
we call the {\bf classifying space} of the bound quiver $(Q,I)$. Several examples are given.

In section $3$, we discuss some homotopy properties of ${\mathcal B}$, and we prove the first main result of this paper.

\subsection*{Theorem \ref{subsec:iso-pi1}}{\em Let $(Q,I)$ be a bound quiver, and ${\mathcal B}= {\mathcal B}(Q,I)$ be its
classifying space. Then  the groups $\pi_1({\mathcal B})$ and $\pi_1(Q,I)$ are isomorphic.}

\smallskip
This generalizes previous similar results obtained independently
for incidence algebras, and for schurian triangular algebras in
\cite{R03} and \cite{Bus02}, respectively. Moreover, this allows
to obtain an adaptation of Van Kampen's theorem to the context of
bound quivers (compare with \cite{R03}).

\medskip In section $4$, we deal with coverings of bound quivers, and of topological spaces. Theorem \ref{subsec:coverings} says that a
covering morphism of bound quivers \mor{(\hat{Q}, \hat{I})}{(Q,I)}{p} induces a covering of topological spaces \mor{{\mathcal B}(\hat{Q},
\hat{I})}{{\mathcal B}(Q,I)}{{\mathcal B}p}. The main result in this section is the following:

\subsection*{Theorem \ref{subsec:galois-coverings}}{\em Let \mor{(\hat{Q},\hat{I})}{(Q,I)}{p}
be a Galois covering given by a group $G$. Then $(\hat{\mathcal B}, {\mathcal B}p)$ is a regular covering of ${\mathcal B}$
with covering automorphism group isomorphic to $G$.}

\medskip In section $5$, we recall the definition of the simplicial (co)homology of an algebra $A=kQ/I$ having a
semi-normed basis. We show that the homotopy relation is strongly related to the vectors of such a basis. Finally, we
compare the simplicial homology and cohomology groups of $A$ to the (cellular) homology and cohomology groups of \BQI \
with coefficients in some abelian group $G$,  ${\rm H}_i({\mathcal B})$ and ${\rm H}^i({\mathcal B}, G)$, respectively.
More precisely, we show the following statement
\subsection*{Corollary \ref{subsec:iso-homology}}{\em Let $A=kQ/I$ be an algebra having a semi-normed basis. Then, for
each $i\geq 0$, there are  isomorphisms of abelian groups
$$\xymatrix@R=2pt{{\rm SH}^i(A,G) \ar[r]^{\sim}   &   {\rm H}^i({\mathcal B}, G), \\
          {\rm SH}_i(A)   \ar[r]^{\sim}   &   {\rm H}_i({\mathcal B}).} $$}

An immediate consequence of this result is that one can think
about the (co)homology theories of \BQI\  as a generalization of
the simplicial (co)homology theories of an algebra $A$, which are
defined only for algebras having semi-normed bases.

\medskip Finally, in section $6$ we focus on the Hochschild cohomology groups of $A$. More precisely, following \cite{GS83,
MdlP99} we compare them with the simplicial cohomology groups of
$A$. Strengthening theorem $3$ in \cite{MdlP99}, we then prove the
following result, which generalizes a theorem due to Gerstenhaber
and Schack \cite{GS83} (see also \cite{C3}).

\subsection*{Theorem \ref{subsec:cohomo-iso}}{\em Let $A = kQ / I$ be a schurian triangular, semi-commutative algebra.
Then, for each $i\geq 0$, there is an isomorphism of abelian
groups $$\xymatrix{{\rm H}^i(\epsilon):{\rm SH}^i(A,k^+) \ar[r]^{\
\ \ \ \ \sim}&{\rm HH}^i(A)}.$$}

It is worth to note that the isomorphisms of Corollary
\ref{subsec:iso-homology} and Theorem \ref{subsec:cohomo-iso} are
induced by isomorphisms of complexes which preserve canonical
cup-products, thus the isomorphism above yields an isomorphism of
graded rings.  We use them to obtain new algebraic-topology
flavored proofs of some known results about the Hochschild
cohomology groups of monomial algebras (\cite{BM98, Hap89}).

%
%
%
%
%
%
%
%
%
%
%
%
%
%
%
%
%
%
%
%
%
%

\section{Preliminaries}\label{sec:Prel}

\subsection{Notation and terminology} Let $Q$ be a finite quiver. We denote by $Q_0$ and $Q_1$ the sets of vertices and
arrows of $Q$, respectively. Given a commutative field $k$, the path algebra $kQ$ is the $k$-vector space with  basis all
the paths of $Q$, including one stationary path $e_x$ for each vertex $x$ of $Q$. Two paths sharing source and target
are said to be {\bf parallel}. The multiplication of two basis elements of $kQ$ is their composition whenever it is
possible, and $0$ otherwise.  Let $F$ be the two-sided ideal of $kQ$ generated by the arrows of $Q$. A two-sided ideal $I$
of $kQ$ is called {\bf admissible} if there exists an integer  $m\geq 2$ such that  $F^m \subseteq I \subseteq F^2$. The
pair $(Q,I)$ is a {\bf bound quiver}. It is well-known that if $A$ is a basic, connected, finite dimensional algebra
over an algebraically closed field $k$, then there exists a unique finite connected quiver $Q$ and a surjective morphism
of $k-$algebras \mor{kQ}{A}{\nu}, which is not unique in general, with $I={\rm Ker}\ \nu$ admissible \cite{BG82}.

Let $(Q,I)$ be a bound quiver, and $A\simeq kQ/I$. It will sometimes be convenient to consider $A$ as a locally bounded
$k-$category, whose object class is $Q_0$, and, for $x,\ y$ in $Q_0$, the morphism set $A(x,y)$ equals the quotient of the
free $k-$module $kQ(x,y)$ with basis the set of paths from $x$ to $y$, modulo the subspace $I(x,y) = I \cap kQ(x,y)$, see
\cite{BG82}. A path $w$ from $x$ to $y$ is said to be a non-zero path if $w \notin I(x,y)$. It is easily seen that
$A(x,y) = e_x Ae_y$.

\subsection{The fundamental group}\label{subsec:fund-group} Given a bound quiver $(Q,I)$, its fundamental group is defined
as follows \cite{MP83}. For $x,\ y$ in $Q_0$, a relation $\rho = \sum_{i=1}^m \lambda_i w_i \in
I(x,y)$ (where $\lambda_i \in k_*$, and $w_i$ are different paths from $x$ to $y$) is said to be {\bf minimal} if
$m\geq 2$, and, for every proper subset $J$ of $\{1, \ldots, n\}$, we have  $\sum_J \lambda_i w_i \notin I(x,y)$.

We define the {\bf homotopy relation} $\sim$ on the set of walks on $(Q,I)$, as the smallest equivalence relation
satisfying

\begin{enumerate}

    \item For each arrow $\alpha$ from $x$ to $y$, one has $\alpha \alpha^{-1} \sim e_x$ and
    $\alpha^{-1} \alpha \sim e_y$.

    \item For each minimal relation $\sum_{i=1}^m \lambda_i w_i$, one has  $ w_i \sim w_j$ for all  $i, j$ in
    $\{1, \ldots, m\}$.

    \item If $u,v, w$ and $w'$ are walks, and $u \sim v$ then $w u w' \sim w v w'$,  whenever these
    compositions are defined.

\end{enumerate}

We denote by $\wt{w}$ the homotopy class of a walk $w$. A closely related notion is that of {\bf natural homotopy}. Two parallel paths
$p$ and $q$ are said to be {\bf naturally homotopic} if $p=q$ or there exists a sequence $p=p_0,\ p_1,\ldots,\ p_s=q$ of parallel paths,
and, for  $i \in \{1,\ldots,\ s\}$, paths $u_i,\ v_i,\ v'_i$ and $w_i$ such that $p_i = u_i v_i w_i$, $p_{i+1} = u_i v'_i w_i$ with $v_i$
and $v'_i$ appearing in the same minimal relation (compare with \cite{AL98}). In that case we write $p\sim_\circ q$ and $\nwt{p}$ will
denote the natural homotopy class of a path $p$. It is easily seen that natural homotopy is the smallest equivalence relation on the set
of paths on $(Q,I)$ satisfying conditions $2)$ and $3)$ (replacing, of course, {\em walks} by {\em paths} in condition $3)$). Thus
$p\sim_\circ q$ implies $p \sim q$ but the converse is not true (see \ref{subsec:def-exemples}, example $1)$). Since the ideal $I$ is
admissible, for every arrow $\alpha$ in $Q$ one has $\nwt{\alpha} = \{\alpha\}$. Moreover, note that if $A=kQ/I$ is schurian, then the
relations $\sim$ and $\sim_\circ$ coincide.

For a fixed vertex $x_0 \in Q_0$,  we denote by $\pi_1(Q,x_0)$ the fundamental group of the underlying graph of $Q$ at the vertex $x_0$.
Let $N(Q,I, x_0)$ be the normal subgroup of $\pi_1(Q, x_0)$ generated by all elements of the form $w^{-1} u^{-1} v w$, where $w$ is a
walk from $x_0$ to $x$, and $u,\ v$ are two homotopic paths from $x$ to $y$. The fundamental group $\pi_1(Q,I)$ is defined to be

$$\pi_1(Q,I) = \pi_1(Q,x_0)  / N(Q,I, x_0).$$

\smallskip Since the quiver $Q$ is connected, this definition is independent of the choice of the base point $x_0$. An
important remark is  that the group defined above depends essentially on the minimal relations, which are given by the
ideal. It is well-known that, for a $k$-algebra $A$, its presentation as a bound quiver algebra is not unique. Thus, the
fundamental group is not an invariant of the algebra (see \cite{A99}). A triangular $k$-algebra $A$ is said to be {\bf
simply connected} if, for every presentation $A  \simeq kQ/I$ of $A$ as a bound quiver algebra, we have $\pi_1(Q,I) = 1$.

However, it has been shown in \cite{BM01} that if an algebra $A$ is {\bf constricted} (that is, if, for every arrow
\mor{x}{y}{\alpha} in $Q_1$, one has ${\rm dim}_k \ A(x, y) = 1$), then the fundamental group is independent of the
presentation.

%
%
%
%
%
%
%
%
%
%
%
%
%
%
%
%
%
%
%
%
%
%

\section{The classifying space \BQI }

\subsection{Background and motivation}\label{subsec:background} To a schurian triangular algebra $A=kQ/I$, is associated a
simplicial complex $|A|$  in the following way \cite{BG83} (see also  \cite{MP84, Bo84}): An $n-$simplex is a sequence
$x_0, x_1,\ldots,x_n$ of $n+1$ different vertices of $Q$ such that for each $j$ with $1\leq j \leq n$, there is a morphism
$f_j$ in $A(x_{j-1}, x_j)$ with $f_n f_{n-1} \cdots f_1 \not=0$. For instance, if $A=A(\Sigma)$ is the incidence algebra of
a poset $(\Sigma, \leq)$, then $|A|$ is the simplicial complex of non empty chains of $\Sigma$ (see \cite{BJ95, GS83}).
Moreover, in this case, this construction is a particular case of that of the classifying space ${\mathcal BC}$ of a small
category ${\mathcal C}$ (see \cite{Se68, Q73}). More generally, let ${\mathcal C}$ be a small category. The space ${\mathcal
B}{\mathcal C}$ is a CW-complex with $0-$cells corresponding to the objects of ${\mathcal C}$, and, for $n\geq 1$, one
$n-$cell for each diagram

$$\xymatrix{X_0 \ar[r]^{f_1} & X_1\ar[r]^{f_2} & \cdots \ar[r]^{f_n}& X_n}$$

\noindent in ${\mathcal C}$ where none of the $f_i$ is an identity map. The corresponding $n-$cell is attached in the obvious
way to any cell of smaller dimension obtained by deleting some $X_i$ and, if $0< i < n$, replacing $f_i$ and $f_{i+1}$
by the composition $f_{i+1}f_i$, whenever this composition is not an identity map. This construction leads to a functor from
the category of small categories to that of CW-complexes.

For example, consider a poset $(\Sigma,\leq)$ as a category whose objects are the elements of $\Sigma$, and, for $x,y\in
\Sigma$ there is a morphism \mor{y}{x}{\rho^y_x} if and only if $x\leq y$ in $\Sigma$, with the obvious
composition. With the above notation, the simplicial complex $|\Sigma|$ is equal to ${\mathcal B}\Sigma$.
This leads us to the following.

\subsection{Definition and Examples} \label{subsec:def-exemples} Recall that, given $n\geq 0$, the standard $n-$simplex is the set
$\Delta^n =\{(t_0, t_1, \ldots, t_n) \in {\mathbb R}^{n+1}| \
\sum_{i=0}^n t_i =1,\ t_i \geq 0\}$. Its $j^{\hbox{\small \it
th}}$ face is  $\partial_j \Delta^n = \{ (t_0, t_1, \ldots, t_n)
\in \Delta^n | \ t_j=0\}$, and, moreover $\partial \Delta^n =
\cup_j  \partial_j \Delta^n$. An open  $n-$cell is an homeomorphic
copy of $\Delta^n \backslash \partial \Delta^n$. Given two
topological spaces $X$, and $Y$, a closet set $A\subset X$, and a
continuous map \mor{A}{Y}{f}, the pushout of $\xymatrix@1{X&A
\ar@{_{(}->}[l]|-{i} \ar[r]|-{f} & Y}$ will be denoted by
$X\coprod_f Y$.

We wish to build a CW-complex by successively attaching $n-$cells to a previously built $(n-1)-$dimensional complex. We begin by giving a
description of the sets ${\mathcal C}_n$ of $n-$cells. Set ${\mathcal C}_0 = Q_0$, ${\mathcal C}_1 =
\{\nwt{\sigma}|\xymatrix@C=12pt{\sigma:x \ar@{~>}[r]& y} \mbox{ is a path in } Q,\ \sigma \notin I, \sigma \not= e_i\}$, and, for $n\geq 2,\ {\mathcal C}_n
= \{(\nwt{\sigma}_1, \ldots, \nwt{\sigma}_n)|\  \sigma_1 \sigma_2\cdots\sigma_n  \mbox{ is a path in } Q,\
\sigma_1\cdots\sigma_n \notin I, \sigma_i \not= e_i \}$. To avoid cumbersome notations, an element $(\nwt{\sigma}_1, \nwt{\sigma}_2, \ldots, \nwt{\sigma}_n)$
of ${\mathcal C}_n$ will be denoted by ${\bm \sigma}^n$, or even ${\bm \sigma}$, if there is no risk of confusion.

Given ${\bm \sigma} = \nwt{\sigma} \in {\mathcal C}_1$, with, say,
$\xymatrix{\sigma:x \ar@{~>}[r]& y}$, define $\partial^1_0({\bm
\sigma}) = y$ and $\partial^1_1({\bm \sigma}) = x.$ More
generally, for $n\geq 2$ and $i\in \{0,\ldots n\}$ define maps
\mor{{\mathcal C}_n}{{\mathcal C}_{n-1}}{\partial^n_i} in the
following way: given ${\bm \sigma} = (\nwt{\sigma}_1,
\nwt{\sigma}_2, \ldots, \nwt{\sigma}_n)$ in ${\mathcal C}_n$, set:
$ \partial^n_0({\bm \sigma})\   =\  (\nwt{\sigma}_2, \ldots ,
\nwt{\sigma}_n),\ldots$,\\   $\partial^n_i({\bm \sigma})\ =\
(\nwt{\sigma}_1, \ldots, \nwt{\sigma_i \sigma_{i+1}},
\ldots,\nwt{\sigma}_n),\ldots,$ and $\partial^n_n({\bm \sigma})\
=\ (\nwt{\sigma}_1, \ldots \nwt{\sigma_{n-1}})$. Again, we shall
write $\partial_i$ instead of $\partial_i^n$. With these notations
we build a CW-complex as follows:

\begin{itemize}

\item $0-$cells: Set ${\mathcal B}_0 = \cup_{x\in Q_0} \Delta^0_x$,

\item $1-$cells: We attach one $1-$cell $\Delta^1_{\bm \sigma}$ for each ${\bm \sigma} \in {\mathcal C}_1$. More precisely, given ${\bm
\sigma} \in {\mathcal C}_1$ define \mor{\partial \Delta^1_{\bm
\sigma}}{{\mathcal B}_0}{f_{\bm \sigma}} by $f_{\bm
\sigma}(\partial_i \Delta^1_{\bm \sigma}) =\Delta^0_{\partial_i^1
{\bm \sigma}}$.  Consider the co-product $f_1 = \coprod_{{\bm
\sigma} \in {\mathcal C}_1} f_{\bm \sigma}$, define $$ {\mathcal
B}_1 = \left( \coprod_{ {\bm \sigma} \in {\mathcal C}_1}
\Delta_{\bm \sigma} ^1\right) \coprod_{f_1} {\mathcal B}_0$$ and
let $p_1$ be the canonical projection \mor{(\coprod_{{\mathcal
C}_1} \Delta^1_{\bm \sigma}) \coprod {\mathcal B}_0}{{\mathcal
B}_1}{p_1}

\item $n-$cells: Assume ${\mathcal B}_{n-1}$ has already been built. Given ${\bm \sigma} \in {\mathcal C}_n$, we have
$\partial_j({\bm \sigma}) \in {\mathcal C}_{n-1}$ for $j\in \{0,\ldots,n\}$. The corresponding $(n-1)-$cells are
$\Delta^{n-1}_{\partial_j {\bm \sigma}}$. Denote by
\mor{\Delta^{n-1}_{\partial_i {\bm \sigma}}}{\coprod_{{\mathcal C}_{n-1}} \Delta^{n-1}_{\bm \tau}}{q_{\partial_i {\bm \sigma}}} the
inclusions, and define \mor{\partial \Delta^n_{\bm \sigma}}{{\mathcal B}_{n-1}}{f_{\bm \sigma}} by
$f_{\bm \sigma}(\partial_j \Delta^n_{\bm \sigma}) = p^{n-1} q_{\partial_j {\bm \sigma}}(\Delta^{n-1}_{\partial_j {\bm \sigma}})$, which is
a continuous function on each $\partial_j \Delta^n_{\bm \sigma}$, thus continuous in $\partial \Delta^n_{\bm \sigma}$. Consider the
co-product
$f_n = \coprod_{{\bm \sigma} \in {\mathcal C}_n} f_{\bm \sigma}$, define
$$ {\mathcal B}_n = \left( \coprod_{ {\bm \sigma} \in {\mathcal C}_n} \Delta_{\bm \sigma} ^n\right) \coprod_{f_n} {\mathcal B}_{n-1}$$
and let $p_n$ be the canonical projection
\mor{(\coprod_{{\mathcal C}_n} \Delta^n_{\bm \sigma}) \coprod {\mathcal B}_{n-1}}{{\mathcal B}_n}{p_n}.

\end{itemize}

Note that since $Q$ is a finite quiver and $I$ is an admissible ideal, there are only finitely many $k\in {\mathbb N}$ such that
${\mathcal C}_k \not= \emptyset$.

\subsection*{Definition}{\em  Let $(Q,I)$ be a bound quiver. The {\em CW-}complex obtained by the preceding construction is the
{\bf classifying space of $(Q,I)$}, and is denoted by ${\mathcal B}(Q,I)$ }

\medskip
A slightly different approach consists in considering homotopy classes of paths, instead of {\em natural} homotopy classes to attach
the cells. The complex obtained in this way will be denoted ${\mathcal B}^\sharp(Q,I)$, and called the {\bf total classifying space of
$(Q,I)$}.

\subsection*{Remarks}
\begin{enumerate}
\item  If the algebra $A=kQ/I$ is {\em almost triangular}, that is $e_x ({\rm rad}\ A) e_y \not= 0$ implies $e_y ({\rm rad}\ A) e_x = 0$ for all vertices
$x, y \in Q_0$ (compare with the definition in \cite{Bus02}), then
the spaces ${\mathcal B}(Q,I)$ and ${\mathcal B}^\sharp(Q,I)$ are
regular CW-complexes.

\item Since for every arrow $\alpha$ of $Q$ we have that $\nwt{\alpha} = \{\alpha\}$, the underlying graph of $Q$ can be considered in a
natural way as a subspace of  the classifying space ${\mathcal B}
(Q,I)$. This is not the case with the total classifying space (see
example 1 below).

\end{enumerate}

\subsection*{Examples}
\begin{enumerate}

    \item Consider the quiver $$\xymatrix{3\ar@/^/[r]^\beta \ar@/_/[r]_\gamma& 2\ar[r]^\alpha &1}$$ bound by $I=<\beta \alpha - \gamma
    \alpha>$. The arrows $\beta$ and $\gamma$ are homotopic, but not naturally homotopic. The spaces
    ${\mathcal B}(Q,I)$ and ${\mathcal B}^\sharp(Q,I)$ look as follows

    \begin{figure}[h]
        \centering
        \resizebox{0.75\linewidth}{!}{
        \input{Exemple1.pstex_t}
        }

    \end{figure}

    \item \label{subsubsec:schurian-case} Let $A=kQ/I$ be a schurian algebra. As noted before, in this case homotopy and natural homotopy
    coincide, thus ${\mathcal B}(Q,I) = {\mathcal B}^\sharp(Q,I)$. For $x,\ y \in Q_0$, there is a $1$-cell joining them if and only if
    there is a non-zero path $w$ from, say, $x$ to $y$. Moreover, if there is another such path $w'$, then, since $A$ is schurian, one has $w
    \sim_\circ w'$. Thus different paths give  the same $1-$cell, and one can identify it with the pair $(x,y)$. In a similar way, given
    an $n-$cell corresponding to $(\nwt{\sigma}_1,\ldots,\nwt{\sigma}_n)$, with $\sigma_i \in A(x_{i-1},x_i)$, one can identify it with
    the sequence of $n+1$ points $x_0,x_1,\ldots ,x_n$ in $Q$. This shows that for schurian algebras $A=kQ/I$, the cellular complex \BQI
    \ is precisely the simplicial complex $|A|$ of \ref{subsec:background}, above (see also \cite{MdlP99}).

    \item Consider the following quiver
    $$\xymatrix@R=12pt{              &6 \ar[d]^\alpha                   &\\
                        &5 \ar[dl]_{\beta_1} \ar[d]_{\beta_2} \ar[dr]^{\beta_3} &\\
2\ar[dr]_{\gamma_1}                     &3 \ar[d]_{\gamma_2}                & 4\ar[dl]^{\gamma_3} \\
                        & 1           &}$$ bound by the ideal generated by $\alpha
                        \beta_1$ and $\sum_{i=1}^3 \beta_i\gamma_i$.
    The cells of \BQI\  are the following:

    \begin{enumerate}

        \item []The $1-$cells are given by
        $\nwt{\alpha},\ \nwt{\beta_i},\ \nwt{\gamma}_i,\ \nwt{\beta_1 \gamma_1} = \nwt{\beta_2 \gamma_2} = \nwt{\beta_3
        \gamma_3},\  \nwt{\alpha \beta_2},$\newline  $\nwt{\alpha \beta_3}$ and
        $\nwt{\alpha \beta_2 \gamma_2} = \nwt{\alpha \beta_3 \gamma_3}$

        \item []The $2-$cells are given by
        $(\nwt{\alpha}, \nwt{\beta_2}),\ (\nwt{\alpha}, \nwt{\beta_3}),\  (\nwt{\alpha},\ \nwt{\beta_2\gamma_2}),
        (\nwt{\alpha},\nwt{\beta_3\gamma_3}),$ \newline
        $(\nwt{\alpha\beta_2},\nwt{\gamma}_2),\ (\nwt{\alpha\beta_3},\nwt{\gamma}_3)$, and
        $(\nwt{\beta_i},\ \nwt{\gamma}_i)$, for $i\in \{1,2,3\}$.

        \item [] The $3-$cells are given by $(\nwt{\alpha},\ \nwt{\beta_2},\ \nwt{\gamma}_2)$ and
        $(\nwt{\alpha},\ \nwt{\beta_3},\ \nwt{\gamma}_3)$.
    \end{enumerate}

    Note that the boundaries of the $2-$cells $(\nwt{\beta}_i,\ \nwt{\gamma}_i)$ are the union of the
    $1-$cells $\nwt{\gamma}_i,\ \nwt{\beta_i},$ and $\nwt{\beta_i \gamma_i}$. Since  $\beta_1 \gamma_1 \sim_\circ
    \beta_2 \gamma_2 \sim_\circ \beta_3 \gamma_3$, the three $2-$cells have a whole $1-$face in common. The
    same argument shows that the two cells of dimension $3$ share a whole face of dimension $2$. The
    geometric realisation of \BQI\ looks as the space shown in the figure below.

    \begin{figure}[h]
\centering

        \includegraphics[width=40mm]{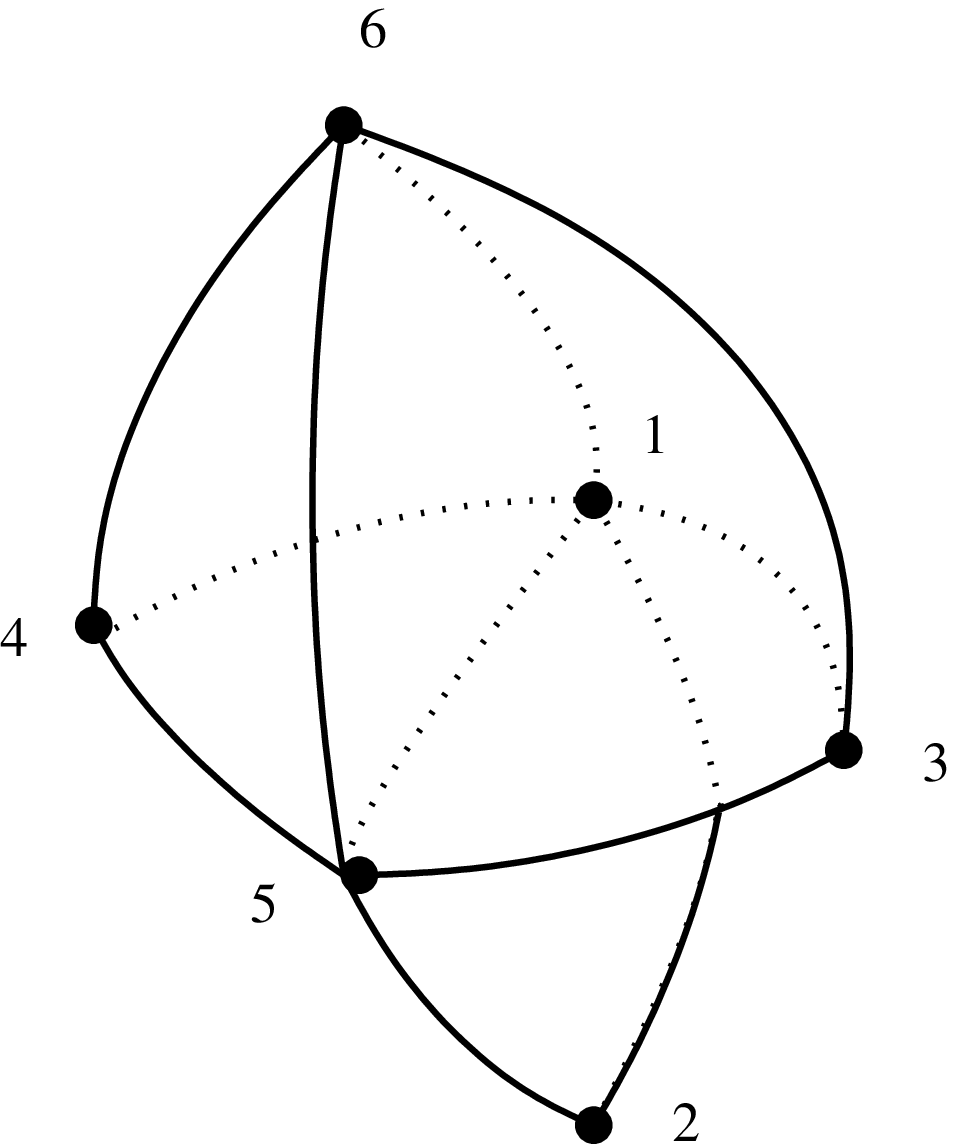}

\end{figure}

    \noindent Note that since $\beta_2 \gamma_2 \sim_\circ \beta_1 \gamma_1$, one has
    $(\nwt{\alpha},\nwt{\beta_2\gamma_2}) = (\nwt{\alpha},\nwt{\beta_1\gamma_1})$, even though the path
    $\alpha \beta_1\gamma_1$ belongs to $I$. However,  $\alpha \beta_1\gamma_1\sim_\circ \alpha \beta_2\gamma_2$ and the
    latter is not zero. Moreover, in this case we have ${\mathcal B}(Q,I)= {\mathcal B}^\sharp(Q,I)$.

    \item \label{subsubsec:radical-carré}Let $Q$ be a quiver, and $I$ be the ideal of $kQ$ generated by paths of length 2. In this case
    there are no minimal relations, so homotopy, and natural homotopy are trivial relations so ${\mathcal B}(Q,I)= {\mathcal
    B}^\sharp(Q,I)$. Each arrow \mor{x}{y}{\alpha} of $Q$ gives rise to a $1-$cell. Moreover, since the only non-zero paths are the
    arrows, these  are the only cells. For the same reason, there are no higher dimensional cells, so the space \BQI \ is homeomorphic to
    $\ol{Q}$, the underlying graph of $Q$.

\end{enumerate}

\subsection*{Remarks}

\begin{enumerate}

    \item  Let ${\mathcal P}(Q)$ denote the path category of $Q$. That is, the object class of ${\mathcal P}(Q)$ is $Q_0$, and for
$x,y \in Q_0$, the morphism set ${\mathcal P}(Q)(x,y)$ is the set of paths from $x$ to $y$ in $Q$. The composition is the obvious
one. Moreover, let ${\mathcal P}(Q,I) ={\mathcal P}(Q) /\sim_\circ$ be the quotient category modulo the natural homotopy relation induced
by $I$. The complex \BQI \ is a subcomplex of the classifying space of ${\mathcal P}(Q,I)$. The n-cells of   ${\mathcal B}({\mathcal
P}(Q,I))$ are in bijection with $n-$tuples $(\nwt{\sigma}_1,\ldots,\nwt{\sigma}_n)$ of composable morphisms, regardless whether their
composition is zero or not. However, if there are no monomial relations in $I$, the complex \BQI \ is  exactly the classifying space
${\mathcal B}({\mathcal P}(Q,I))$. The same applies for  ${\mathcal B}^\sharp(Q,I)$ with respect to the classifying space of the category
${\mathcal P}(Q)/\sim$.

    \item  Consider the quiver
$$\xymatrix@C=12pt@R=7pt{
            &           &           &  5 \ar[drr]       &   &\\
            &2 \ar[drr]\ar[urr]&           &           &   &8\\
            &           &4\ar[drr]\ar[ruu]  &  6 \ar[urr]  &   &\\
1\ar[uur]\ar[urr]\ar[drr]&          &           &           &7\ar[uur]&\\
            &           &3\ar[urr]\ar[uur]&            &   &}$$

\noindent bound by all the commutativity relations and all the paths of length 3. The space \BQI \ is homeomorphic to the
$3-$dimensional sphere $S^2$. The commutativity relations tell how to "glue" the $2-$dimensional cells. The existence of
monomial relations implies that there are no $3-$cells to "fill the hole". On the other hand, the space   ${\mathcal
B}({\mathcal P}(Q,I))$ is homeomorphic to the ball $B^3$.

Let $(Q,I')$ be the same quiver bound by the commutativity relations (and only these). In particular, the natural homotopy
relations induced by $I$ and $I'$ are the same. However,  the space ${\mathcal B}(Q,I')$ is homeomorphic to $B^3$, hence
does not have the same homotopy type as \BQI. This shows that monomial relations play an important rôle in the construction
of \BQI, even though they are  taken into account to define the neither the natural homotopy nor the homotopy relations.

    \item  The fact that ${\mathcal B}(Q,I)$ is not really the classifying space of a category implies that this
construction is not functorial, as the following example shows: Consider the quiver $Q$
$$\xymatrix@R=1pt@C=25pt{  & 3 \ar[dr]^{\beta_1}&\\
4 \ar[dr]_{\alpha_1} \ar[ur]^{\beta_1}&   &1\\
                    & 2\ar[ur]_{\alpha_2}& }$$

\noindent and let $I= <\alpha_1 \alpha_2 + \beta_1 \beta_2>,\  I'= <\alpha_1 \alpha_2, \beta_1 \beta_2>$.
The identity map on $Q$ yields a bound quiver morphism \mor{(Q,I)}{(Q,I')}{f}. However the induced map \mor{(Q,I')}{(Q,I)}{g} is not
a bound quiver morphism since $\alpha_1 \alpha_2 \in I'$, but $\alpha_1 \alpha_2 \notin I$. In this case, the space ${\mathcal B}(Q,I)$ is
homeomorphic to $B^2=\{x \in{\mathbb R}^2|\ ||x|| \leq1\}$. On the other hand, ${\mathcal B}(Q,I')$ is homeomorphic to $S^1$. Thus, there
is no non homotopically trivial map from ${\mathcal B}(Q,I)$ to ${\mathcal B}(Q,I')$

\end{enumerate}

%
%
%
%
%
%
%
%
%
%
%
%
%
%
%
%
%
%
%
%
%
%

\section{Homotopy}

\subsection{Proposition} \label{subsec:retract}{\em Let $(Q,I)$ be a bound quiver with $I$ a monomial ideal such that $=kQ/I$ is almost
triangular. Then the underlying graph  $\overline{Q}$ is a deformation retract of \BQI.}

\pf Since $I$ is monomial, the natural homotopy and the homotopy
relations are trivial. Thus, the $n-$cells are given by $n-$tuples
$(\sigma_1,\ldots,\sigma_n)$ of paths such that $\sigma_1
\sigma_2\cdots\sigma_n \notin I$. If $I = F^2$, there is nothing
to prove (see example \ref{subsubsec:radical-carré}, above). If
this is not the case, then there are non-zero paths of length
greater or equal than $2$. Let $\alpha_1 \cdots \alpha_n$ be a
maximal non-zero path. It gives a $n-$cell ${\bm \alpha}$ which is
maximal. For $i\in\{0,\ldots,n-1\}$, let $x_i$ be the source of
$\alpha_{i+1}$, and $x_n$ be the target of $\alpha_n$. Since $A$
is almost triangular, $x_i \not= x_j$ whenever $i\not=j$, thus we
may consider the $n-$cell ${\bm \alpha}$ as the standard
$n-$simplex $[x_0,\ldots x_n]$. Let $y=\frac{1}{2}(x_0 + x_n)$,
thus $[x_0,\ldots x_n] = [x_0,\ldots,x_{n-1},y] \cup
[y,x_1,...,x_n]$.

\begin{figure}[h]
        \centering
        \resizebox{0.85\linewidth}{!}{
        \input{prop31.pstex_t}
        }
\end{figure}

For $t\in [0,1]$ define $r(t) = (1-t)y + tx_1$, that is the path
joining $y$ to $x_1$. Using $r$ one can "crush" the simplex
$[x_0,\ldots x_n]$ onto its faces $[x_0,\ldots,x_{n-1}]$ and
$[x_1,\ldots,x_n]$, which, by the maximality assumption, are the
only faces in the complex that are not free (use the barycentric
coordinates, and the division of the simplex given below). The
space obtained is \BQI, in which we have crushed the cell ${\bm
\alpha}$ onto the $(n-1)-$cells $\partial_0 ({\bm \alpha})$ and
$\partial_n({\bm \alpha})$. It is easily seen that this space is
precisely ${\mathcal B}(Q,I_1)$ where $I_1 = I + <\alpha_1 \cdots
\alpha_n>$, which is again a monomial ideal. If $I_1 = F^2$, we
have finished, otherwise we repeat the process above with another
maximal dimensional cell of ${\mathcal B}(Q,I_1)$. This must end
in a finite number of steps, since $I$ is admissible.\cqfd

A first immediate consequence is that if $I$ is a monomial ideal, then the fundamental group $\pi_1({\mathcal B})$ of the
topological space ${\mathcal B}$,  is isomorphic to the fundamental group $\pi_1(Q,I)$ of the bound quiver  $(Q,I)$. In
fact, this is a particular case of a much more general situation.

Given a cell complex $X$, its fundamental group $\pi_1(X)$ can be
described in the following convenient way (see \cite{MASS89}, for
instance).  Fix a $0-$cell, $x_0$, and a maximal tree $M$, that
is, a subcomplex of dimension smaller or equal than $1$, which is
acyclic and maximal with respect to this property. For every
$2-$cell $e^2_{\lambda}$ of $X$, let $\alpha_e$ be a path class
which starts at some fixed $0-$cell $x$ in its boundary, $\partial
e^2_{\lambda}$, and goes around $\partial e^2_{\lambda}$ exactly
once. Moreover, let $\beta_x$ be the unique path class in $X$
which goes from $x_0$ to $x$ along the tree $M$. Finally set  $
\gamma_e = \beta_x \alpha_e \beta_x^{-1}$ (compare with the
definition of a {\bf parade data} in \cite{FGM00}).

Let $G$ be the free group on the set of $1-$cells of $X$, and $N$ be the normal subgroup of $G$ generated by the following
elements

\begin{enumerate}

    \item The cells of $M$,

    \item The elements $\gamma_e$, as constructed above.

\end{enumerate}

With the above notation, reformulating theorem 2.1, p 213 of \cite{MASS89}, one has  $\pi_1(X) \simeq G/N$.

Moreover, since we are interested in complexes of the form ${\mathcal B} =  {\mathcal B}(Q,I)$, and in this case all the
$2-$cells are of the form $(\nwt{\sigma}_1, \nwt{\sigma}_2)$, the boundary of such a cell being $\nwt{\sigma}_2, \
\nwt{\sigma_1 \sigma_2}, \  \nwt{\sigma}_1$, we can improve the preceding presentation for $\pi_1({\mathcal B}).$

\subsection{Lemma}{\em Let $(Q,I)$ be a bound quiver and $T$ be a maximal tree in $Q$.
Then $\pi_1({\mathcal B}) \simeq F/K$, where $F$ is the free group with basis the set of arrows of $Q$, and $K$ is the
normal subgroup of $F$ generated by the elements of the following two types:

\begin{enumerate}

    \item $\alpha$, for every arrow $\alpha$ in $T$

    \item $(\alpha_1 \alpha_2 \cdots \alpha_r)(\beta_1 \beta_2 \cdots \beta_s)^{-1}$ whenever  $\alpha_1
    \alpha_2 \cdots \alpha_r$  and $\beta_1 \beta_2 \cdots \beta_s$ are two paths appearing in the same minimal relation.

\end{enumerate} }

\pf Let $T$ be a maximal tree in the quiver $Q$. In particular $T$ is a set of arrows of $Q$. It follows from the construction of \BQI \
that an arrow \mor{x}{y}{\alpha} in $T$ gives rise to a $1-$cell $\nwt{\alpha}$ in \BQI. The set of $1-$cells obtained from the arrows of
$T$ forms a maximal tree $M$ in \BQI. We work with respect to these maximal trees. Moreover, let $G$ and $N$ be as before.

Consider the map \mor{F}{G/N}{\phi} defined by \appl{\alpha}{\tilde{\alpha}N}. It is straightforward to see that
$ K \subseteq {\rm Ker}\ \phi$, so we obtain a group homomorphism \mor{F/K}{G/N}{\Phi} defined by \appl{\alpha K}{\tilde{\alpha}N}

We show that $\Phi$ is an isomorphism by constructing its inverse.
If $\nwt{w}$ is a $1-$ cell in \BQI, then there are arrows
$\alpha_1, \alpha_2,\ldots,\alpha_r$ in $Q$ such that $\alpha_1
\alpha_2\cdots\alpha_r \sim_\circ w$, and $\alpha_1
\alpha_2\cdots\alpha_r \notin I$. Define the map
\mor{G}{F/K}{\psi} by \appl{\tilde{w}}{\alpha_1
\alpha_2\cdots\alpha_r K}. It is not hard to see that this is a
well-defined map, and, moreover, that $N \subseteq {\rm Ker}\
\psi$. Thus, we obtain a group homomorphism \mor{G/N}{F/K}{\Psi}
defined by \appl{\nwt{w}N}{\alpha_1 \ldots \alpha_r K}. Finally,
it is straightforward to check that $\Phi$ and $\Psi$ are mutually
inverse.\cqfd

\subsection{Theorem}\label{subsec:iso-pi1} {\em Let $(Q,I)$ be a bound quiver, with $Q$ triangular and ${\mathcal B} =
{\mathcal B}(Q,I)$. Then the groups $\pi_1({\mathcal B})$ and $\pi_1(Q,I)$ are isomorphic.}

\pf This is an easy consequence of the preceding lemma, and the description of $\pi_1(Q,I)$ given in section
\ref{subsec:fund-group}.\cqfd

\subsection*{Remark} A similar argument applies to the fundamental group of ${\mathcal B}^\sharp(Q,I)$. Thus we obtain
$\pi_1({\mathcal B}^\sharp) \simeq \pi_1(Q,I)$.

\subsection{Corollary}{\em Let $A=kQ/I$ be a triangular algebra. Then $A$ is simply connected if and only if
for every presentation $(Q,I)$,  \BQI \ or, equivalently ${\mathcal B}^\sharp (Q,I)$,  is a simply connected topological space.} \cqfd

Given a bound quiver $(Q,I)$, and a full convex subquiver $Q^i$, let $I^i$ denote the ideal $I \cap kQ^i$ of $kQ^i$, that is, the
restriction of $I$ to $Q^i$. With these notations we have the following result (see also \cite{R03}).

\subsection{Corollary} {\em Let $(Q,I)$ be a bound quiver, $Q^1$ and $Q^2$ be two full convex subquivers of $Q$ such that every
non-zero path of $Q$ lies entirely in either $Q^1$ or $Q^2$, and $Q^0 = Q^1 \cap Q^2$ is connected. Then,
$\pi_1(Q,I)$ is the pushout of the diagram
$$\xymatrix@C=25pt@R=20pt{
\pi_1(Q^1,I^1)& \pi_1(Q^0,I^0) \ar[l] \ar[r] & \pi_1(Q^2, I^2)}$$ where the arrows are the maps induced by the inclusions.}

\pf It follows from the hypotheses made on $Q^1$ and $Q^2$, that ${\mathcal B}(Q^1, I^1) \cup {\mathcal B}(Q^2, I^2) = {\mathcal
B}(Q,I)$, and that ${\mathcal B}(Q^1, I^1) \cap {\mathcal B}(Q^2, I^2)$ is connected. The result then follows from theorem
\ref{subsec:iso-pi1}, and Van Kampen theorem for topological spaces (see \cite{ROT88}, for instance). \cqfd

\subsection*{Example} Consider the quiver $Q$

$$\xymatrix{& 1&\\
3\ar@/^/[r]^{\alpha_1} \ar@/_/[r]_{\beta_1}&    2 \ar@/^/[u]^{\alpha_2} \ar@/_/[u]_{\beta_2} \ar[r] \ar[dr]  &4\\
                            &6 \ar[ur] \ar[r]                       &5}$$
bound by $I=<\alpha_1 \alpha_2 - \beta_1 \beta_2, \alpha_1 \beta_2 - \beta_1 \alpha_2>$, and let $Q^1$ and $Q^2$ be the full
subquivers of $Q$ generated be the vertices $2,3,4,5,6$, and $1,2,3$, respectively. An easy computation gives
$\pi_1(Q^1, I^1) = {\mathbb Z} \ast {\mathbb Z},\  \pi_1(Q^2, I^2) = {\mathbb Z}_2$, and
$\pi_1(Q^0, I^0) = {\mathbb Z}$, and this yields $\pi_1(Q, I) = {\mathbb Z} \ast {\mathbb Z}_2$.

%
%
%
%
%
%
%
%
%
%
%
%
%
%
%
%
%
%
%
%
%
%

\section{Coverings}

As we saw in the example before proposition \ref{subsec:retract}, the construction of \BQI \ does not lead to a
functor from the category of bound quivers with bound quivers morphisms to the category of CW-complexes with cellular
maps. However, as we now see, if we consider the category of bound quivers with covering morphisms, then ${\mathcal B}$
gives rise to a covering morphism of topological spaces.

\medskip

We refer the reader to \cite{BG82, G80, P86}, for references about coverings and their uses in representation
theory of algebras.

Given a bound quiver $Q$ and a vertex $x \in Q_0$, let, as usual, $x^+$ (or $x^-$) be the set of all arrows leaving (or
entering, respectively) the vertex  $x$. A bound quiver morphism \mor{(\wh{Q},\wh{I})}{(Q,I)}{p}  induces in an obvious
way a $k-$linear map \mor{k\wh{Q}}{kQ}{p}. Recall from \cite{P86}, for instance,  that a bound quiver morphism
\mor{(\wh{Q},\wh{I})}{(Q,I)}{p} is called a {\bf covering morphism} if

\begin{enumerate}

    \item $p^{-1}(x) \not= \emptyset $ for every $x\in Q_0$

    \item For every $x \in Q_0$, and  $\hat{x} \in p^{-1}(x), $ the map $p$ induces bijections
    $\xymatrix@1@C=15pt{\hat{x}^+ \ar[r]&x^+}$ and $\xymatrix@1@C=15pt{\hat{x}^- \ar[r]&x^-}$

    \item For every $x,y\ \in Q_0,$ every relation $\rho \in I(x,y) $ and every $\hat{x} \in p^{-1}(x)$ there exists
    $\hat{y} \in p^{-1}(y)$ and $\hat{\rho} \in \hat{I}(\hat{x},\hat{y})$ such that $p(\hat{\rho}) = \rho$.

\end{enumerate}

Conditions $2)$ and $3)$  ensure that covering morphisms behave
well with respect to homotopy relations. We have the following
straightforward lemmata whose easy proofs will be omitted.

\subsection{Lemma}{\em
Let \mor{(\wh{Q},\wh{I})}{(Q,I)}{p} be a covering morphism, $x$ be a vertex of $Q$, and $w_1, \ w_2$ be two paths with
source $x$. Moreover, let $\hat{x} \in p^{-1}(x)$, and $\wh{w}_1,\ \wh{w}_2$ be two paths with source $\hat{x}$ such that
$p(\wh{w}_i)= w_i$, for $i = 1,2$. Then $w_1 \sim_\circ w_2$ if and only if $\wh{w}_1 \sim_\circ \wh{w}_2$.}\cqfd

\subsection{Lemma}{\em  Let \mor{(\wh{Q},\wh{I})}{(Q,I)}{p} be a covering morphism. For $x_0 \in Q_0$, and every $\hat{x}_0 \in
p^{-1}(x_0)$, there is a bijective correspondence between the set of $n-$cells of ${\mathcal B} = {\mathcal B}(Q,I)$ having
$x_0$ as boundary point, and the set of $n-$cells of $\wh{\mathcal B} = {\mathcal B}(\wh{Q},\wh{I})$ having $\hat{x}_0$ as
boundary point.}\cqfd

In light of the preceding result, we can define \mor{\hat{\mathcal
B}}{{\mathcal B}}{{\mathcal B}p} as the map which maps
homeomorphically an $n-$cell $(\nwt{s}_1,\ldots \nwt{s}_n)$ onto
the cell $(\nwt{p(s_1)}, \ldots \nwt{p(s_n)})$.

Recall from \cite{ROT88}, for instance, that if $X$ is a topological space, then a covering space of $X$ is a pair $(\hat{X},p)$ where

\begin{enumerate}

    \item $\hat{X}$ is an arcwise connected topological space

    \item \mor{\hat{X}}{X}{p} is a continuous map

    \item Each $x \in X$ has an open neighborhood $U_x$ such that $p^{-1}(U_x) = \bigcup_{i\in I} \hat{U}_i$, with
    $\hat{U}_i$ disjoint open sets, and \mor{\hat{U}_i}{U_x}{p|_{\hat{U}}} an homeomorphism, for every $i\in I$.

\end{enumerate}

This yields the following result.

\subsection{Theorem}\label{subsec:coverings}{\em Let \mor{(\wh{Q},\wh{I})}{(Q,I)}{p} be a
covering morphism of bound quivers, with $Q$, $\wh{Q}$ connected.
Then $(\hat{\mathcal B}, {\mathcal B}p)$ is a covering space of
\BQI}

\pf Since $\hat{Q}$ is a connected quiver, $\wh{\mathcal B}$ is an
arcwise connected space. Moreover, it follows from its definition
that ${\mathcal B}p$ is a continuous map. Thus, there only remains
to prove that the third condition of the definition is satisfied.
But this follows from the fact that the open cells of
$\wh{\mathcal B}$ are disjoint open sets, and from the definition
of ${\mathcal B}p$, whose restriction the each such cell is an
homeomorphism.\cqfd

Among all the covers of a bound quiver $(Q,I)$, there is one of
particular interest. The {\bf universal cover} of $(Q,I)$ is a
cover map \mor{(\hat{Q}, \hat{I})}{(Q,I)}{p} such that for any
other cover \mor{(\ol{Q}, \ol{I})}{(Q,I)}{p'} there exists a
covering map \mor{(\hat{Q}, \hat{I})}{(\ol{Q}, \ol{I})}{\pi}
satisfying $p=p'\pi$ (see \cite{MP83}).

\subsection{Corollary}{\em  If \mor{(\hat{Q},\hat{I})}{(Q,I)}{p} is the universal cover of $(Q,I)$, then
$(\wh{\mathcal B}, {\mathcal B}p)$ is the universal cover of ${\mathcal B}$.} \cqfd

Recall from \cite{P86}, for instance, that a cover of bound quivers \mor{(\hat{Q},\hat{I})}{(Q,I)}{p} is said to be a
{\bf Galois cover} defined by the action of a group $G$ of automorphisms of $(\hat{Q},\hat{I})$ if

\begin{enumerate}

    \item [$(4)$] $ pg = p$ for all $g\in G$

    \item [$(5)$] $ p^{-1}(x) = G p^{-1}(x)$ and $p^{-1}(\alpha) = G p^{-1}(\alpha)$, for all vertices $x$ in $Q_0$ and arrows
    $\alpha \in Q_1.$

    \item [$(6)$] $ G$ acts without fixed points on $\hat{Q}$.

\end{enumerate}

Moreover, in this situation there exists a normal subgroup $H$ of $\pi_1(Q,I)$ such that $\pi_1(\hat{Q},\hat{I}) \simeq H$ and
$\pi_1(Q,I) / H \simeq G$ (see \cite{P86}).

As before, an automorphism $g$ of $(\hat{Q},\hat{I})$ induces a
map from the set of paths of $(\hat{Q},\hat{I})$ to the set of
paths of $(Q,I)$. This allows to define a cellular map
\mor{\hat{\mathcal B}}{\mathcal B}{{\mathcal B}g} as the
continuous function that maps homeomorphically the cell
$(\nwt{\sigma}_1, \ldots, \nwt{\sigma}_n)$ onto the cell
$(\nwt{g\sigma_1}, \ldots, \nwt{g\sigma_n})$. The fact that
${\mathcal B}g$ is well-defined follows from the fact that $g$ is
invertible. Moreover, it is straightforward to check that
${\mathcal B}g {\mathcal B}p$ = ${\mathcal B}p$. An important
remark is that the restriction of ${\mathcal B}g$ to the $0-$cells
of $\hat{\mathcal B}$ is precisely $g$. An immediate consequence
of this is that if $g_1$ and $g_2$ are automorphisms of
$(\hat{Q},\hat{I})$ with $g_1 \not=g_2$, then ${\mathcal B}g_1
\not= {\mathcal B}g_2$.

On the other hand, given a covering space $(\hat{X}, p)$ of a
topological space $X$, we denote  by $p_*$ the group homomorphism
\mor{\pi_1(\hat{X})}{\pi_1(X)}{\pi_1(p)}. In particular, $p_*$ is
always a monomorphism (see \cite{ROT88}). If ${\rm Im}\ p_*$ is a
normal subgroup of $\pi_1(X)$ then the covering $(\hat{X},p)$ is
said to be {\bf regular}. The set of all homeomorphisms
\mor{\hat{X}}{\hat{X}}{\phi} such $p\phi = p$ is a group, which is
called the group of covering automorphisms of $(\hat{X},p)$, and
is denoted by ${\rm Cov}(\hat{X}/X)$. It is well-known that a
covering $(\hat{X},p)$ is regular if and only if ${\rm
Cov}(\hat{X}/X)$ acts transitively on $p^{-1}(x_0)$, the fiber
over the base point (see \cite{ROT88}, for instance).

Note that the set $\{{\mathcal B}g| \ g\in G\}$ is a subgroup of ${\rm Cov}(\hat{\mathcal B}/{\mathcal B})$, which is
isomorphic to $G$. In fact, we have the following stronger result.

\subsection{Theorem}\label{subsec:galois-coverings}{\em Let \mor{(\hat{Q},\hat{I})}{(Q,I)}{p}
be a Galois covering given by a group $G$. Then $(\hat{\mathcal B}, {\mathcal B}p)$ is a regular covering of ${\mathcal B}$
and ${\rm Cov}(\hat{\mathcal B}/{\mathcal B}) \simeq G$.}

\pf First of all, fix a vertex $x_0 \in Q_0$, which is also a $0-$cell of  ${\mathcal B}$. These will be the base points
with respect the fundamental groups that will be considered.

In order to show that $(\hat{\mathcal B}, {\mathcal B}p)$ is a
regular covering of ${\mathcal B}$, it is enough to show that
${\rm Cov}(\hat{\mathcal B}/{\mathcal B})$ acts transitively on
the fiber $p^{-1}(x_0)$. This follows immediately from condition
$5$ in the definition of a Galois covering, and the fact that $G$
is isomorphic to a subgroup of ${\rm Cov}(\hat{\mathcal
B}/{\mathcal B})$.

On the other hand

$${\rm Cov}(\hat{\mathcal B}/{\mathcal B}) \simeq \frac{\pi_1({\mathcal B})}{({\mathcal B}p)_* \pi_1(\hat{\mathcal B})}
\simeq \frac{\pi_1(Q,I)}{ \pi_1(\hat{Q}, \hat{I})} \simeq G $$
\noindent where the first isomorphism is given by Corollary 10.28, p. 294 in \cite{ROT88}\cqfd

\subsection*{Remark} Theorem 10.19, p.290, in \cite{ROT88} states that a nontrivial covering automorphism
$h \in {\rm Cov}(\hat{X}/X)$ acts without fixed points on $\hat{X}$. Thus, condition $6$ in the definition of Galois
coverings is redundant.

\subsection*{Example} Consider the quiver $Q$
$$\xymatrix{ 3\ar@/^/[r]^{\alpha_1} \ar@/_/[r]_{\beta_1}& 2\ar@/^/[r]^{\alpha_2} \ar@/_/[r]_{\beta_2}&1}$$

\noindent bound by the ideal $I=<\alpha_1 \alpha_2 - \beta_1 \beta_2, \alpha_1 \beta_2 - \beta_1 \alpha_2>$. The space \BQI
\ is homeomorphic to the real projective plane ${\mathbb R}P^2$. On the other hand, the universal cover of $(Q,I)$ is given by
the quiver $Q'$

$$\xymatrix{ x_3 \ar[r] \ar[dr]  & x_2\ar[r] \ar[dr]& x_1\\
         y_3 \ar[r]\ar[ur]  & y_2\ar[r] \ar[ur]& y_1}$$

\noindent bound by the ideal $I'$ generated by all the possible commutativity relations. The covering map \mor{Q'}{Q}{p} is given by
$p(x_i)= p(y_i) = i$, for $i\in \{1,2,3\}$. The space ${\mathcal B}(\wh{Q},\wh{I})$ is homeomorphic to the sphere
$S^2$, and the map ${\mathcal B}p$ identifies antipodal points. It is well-known that the pair $(S^2, {\mathcal B}p)$ is a
covering space of ${\mathbb R}P^2$.

\subsection*{Remark} Again, with the obvious changes in lemmata $4.1$, $4.2$ and theorem $4.3$ we obtain a similar result for
the total classifying space ${\mathcal B}^\sharp(Q,I)$, and its
covering spaces.

%
%
%
%
%
%
%
%
%
%
%
%
%
%
%
%
%
%
%

\section{(Co)Homologies}

\subsection{(Co)Homology of \BQI} Let $({\rm C}_{\bullet}({\mathcal B}), \delta)$ be the complex defined by letting ${\rm C}_i({\mathcal
B})$ be the free abelian group with basis the set of $i-$cells, ${\mathcal C}_i$, and by defining \mor{{\rm C}_n({\mathcal B})}{{\rm
C}_{n-1}({\mathcal B})}{\delta_n} on the basis elements by the rule $\delta_n({\bm \sigma})   = \sum_{i=0}^{n+1}(-1)^i \partial^n_i ({\bm
\sigma})$.

With these notations, the homology groups  ${\rm H}_i({\mathcal
B})$ of ${\mathcal B}$ are the homology groups of $({\rm
C}_{\bullet}({\mathcal B}), \delta)$. Moreover, if $G$ is an
abelian group, then the cohomology groups ${\rm H}^i({\mathcal B},
G)$ of ${\mathcal B}$ with coefficients in $G$, are the cohomology
groups of ${\rm Hom}_{\mathbb Z}({\rm C}_{\bullet}({\mathcal
B}),G)$. In an analogous way, we define the corresponding complex
$({\rm C}_{\bullet}({\mathcal B}^\sharp), \delta)$, and obtain the
(co)homology groups of ${\mathcal B}^\sharp$.

\subsection{Simplicial (co)homology of algebras}

On the other hand, recall that the simplicial homology of a
schurian algebra was defined in \cite{BG83} (see also
\cite{MdlP99, MP84} for generalization to the non-schurian case)
in the following way. For every pair of vertices $i,j$ of $Q$, let
$B_{ij}$ be a $k-$basis of $e_i A e_j$. The set $B = \cup_{i,j}
B_{ij}$ is said to be a {\bf semi-normed basis} if:

\begin{enumerate}

    \item $e_i \in B_{ii}$, for every vertex $i$.

    \item $\ol{\alpha} \in B_{ij}$, for every arrow \mor{i}{j}{\alpha}

    \item For $\sigma, \sigma'$ in $B$, either $\sigma \sigma'= 0$ or there exist $\lambda_{\sigma \sigma'} \in k$,
    and $b(\sigma, \sigma') \in B$ such that $\sigma \sigma' = \lambda_{\sigma \sigma'} b(\sigma, \sigma')$
\end{enumerate}

Define the chain complex $({\rm SC}_{\bullet}(A), d)$ in the following way: ${\rm SC}_0(A)$, and ${\rm SC}_1(A)$ are the free  abelian
groups with basis $Q_0$, and $B$, respectively. For $n\geq 2$, let ${\rm SC}_n(A)$ be the free abelian group with  basis the set of
$n-$tuples $(\sigma_1,\ldots,\sigma_n)$ of $B^n$ such that $\sigma_1 \sigma_2 \cdots \sigma_n \not= 0$, and $\sigma_i \not= e_j$, for all
$i$, for all $j\in Q_0$. The differential is induced by the rules $d_1(\sigma) = y-x$ when $\sigma \in B_{xy}$, and, for $n\geq 1$,

\begin{eqnarray*}
d_n(\sigma_1 \sigma_2 \cdots \sigma_n) &=& (\sigma_2,\ldots,\sigma_n)\\
            & + &\sum_{j=1}^{n-1} (-1)^j (\sigma_1,\ldots, b(\sigma_j,\sigma_{j+1}),\ldots,\sigma_n)\\
            & +& (-1)^n (\sigma_1,\ldots,\sigma_{n-1})
\end{eqnarray*}

\noindent where  $b(\sigma_j, \sigma_{j+1})$ is as in condition $3$, above. The $n^{th}$ homology group ${\rm SH}_n(A)$ of   $({\rm
SC}_{\bullet}(A, d))$ is called the {\bf $n^{th}$ simplicial homology group of $A$} with respect to the basis $B$. For an  abelian group
$M$, the $n^{th}$ cohomology group ${\rm SH}^n(A, M)$, of the cochain complex ${\rm Hom}_{\mathbb Z}({\rm SC}_{\bullet}(A),M)$, {\bf
$n^{th}$ simplicial cohomology group of $A$}, with coefficients in $M$, with respect to the basis $B$. Moreover, by $d^n$, we denote the
induced differential ${\rm Hom}_{\mathbb Z}(d_n, M)$. \subsection*{Remarks}

\begin{enumerate}

    \item Not every algebra of the form $A=kQ/I$ admits a semi-normed basis. However, it is shown in \cite{MdlP99}
    that if the universal cover of $(Q,I)$ is schurian, then $A$ admits a semi-normed basis. In particular all schurian
    algebras admit semi-normed basis.

    \item On the other hand, $A$ has a semi-normed basis $B$, the groups ${\rm SC}_i$ depend essentially on the way
    $B$ is related to $I$. Hence, as for fundamental groups, different presentations of the algebra  may lead to
    different simplicial homology groups, as the following well-known example shows: Consider the quiver $Q$
    $$\xymatrix{3\ar@/^/[r]^\beta \ar@/_/[r]_\gamma& 2\ar[r]^\alpha &1}$$ Let $I_1$ be the ideal generated by $\beta
    \alpha$, and $A_{\nu_1} = kQ/ I_1$.  Associated to this presentation,  there is a semi-normed basis $\{e_1, e_2,
    e_3, \ol{\alpha},\ol{\beta},\ol{\gamma}, \ol{\gamma\alpha} \}$, and with respect to this basis,  ${\rm
    SC}_2(A_{\nu_1})$ has as basis $\{(\gamma, \alpha)\}$, and this leads to  ${\rm SH}_1(A_{\nu_1}) = {\mathbb Z}$.

    On the other hand, let $I_2$ be the ideal generated by $(\beta - \gamma)\alpha$, and $A_{\nu_2}= kQ/I_2$. It is
    easy to see that $A_{\nu_1} \simeq A_{\nu_2}$. With this presentation, one has the same semi-normed basis, however ${\rm
    SC}_2(A_{\nu_2})$ has a basis $\{(\ol{\gamma}, \ol{\alpha}), (\ol{\beta},\ol{\alpha})\}$, and this leads to  ${\rm
    SH}_1(A_{\nu_2}) = 0$.

    \item As noted in \cite{MdlP99}, in case $A$ is schurian one can identify a basis element
    $(\sigma_1,\ldots,\sigma_n)$ of ${\rm SC}_n(A)$ with the $(n+1)-$tuple $(x_0, x_1,\ldots, x_n)$ of vertices of $Q$,
    where $\sigma_i \in e_{x_{i-1}}Ae_{x_i}$, for $1\leq i \leq n$ (compare with example \ref{subsubsec:schurian-case} in
    section \ref{subsec:def-exemples}). Thus, for schurian algebras the simplicial homology groups are independent of
    the semi-normed basis with respect to which they are computed. Moreover, it is straightforward that in this case one
    has, for every $i\geq 0$, ${\rm H}_i({\mathcal B}) \simeq {\rm SH}_i(A)$,
    and ${\rm H}^i({\mathcal B}, G) \simeq {\rm SH}^i(A,G)$ for every abelian group $G$.

\end{enumerate}

The last remark is in fact a particular case of a more general result. Before stating and proving it we need the following
technical lemma.

\subsection{Lemma} \label{subsec:lemma_bases}{\em Let $A=kQ/I$ be an algebra having a semi-normed basis $B$, then:

    \begin{enumerate}

    \item [$a)$] For every non-zero path $w$ there exist a unique basis element $b(w)$ and a scalar $\lambda$ such
    that $\ol{w} = \lambda b(w)$.

    \item [$b)$] For every basis element $v \in B$, there exist a non-zero path $p(v)$ and a scalar $\mu$ such that
    $v = \mu \ol{p(v)}$. Moreover, if $p(v)$ and $p'(v)$ are two such paths, then they are naturally homotopic.

    \item [$c)$] If $w_1, w_2$ are two non-zero paths with $w_1 \sim_\circ w_2$, then there
    is a scalar $\lambda \in k$ such that $\ol{w}_2 = \lambda \ol{w}_1$

    \end{enumerate}}
\pf
\begin{enumerate}

\item [$a)$] This  follows from an easy induction on the length  $l(w)$ of the path $w$.

\item [$b)$]  Let $\sigma \in B_{xy}$, and  $w_1,\  w_2,\ldots,\ w_r$ be paths from $x$ to $y$ such that $\{\ol{w_i}|\ 1\leq i \leq r \}$
is a basis of $e_x A e_y$. It follows from $a)$ that there exists
$\lambda_i \in k\backslash \{0\}$ and $b(w_i) \in B_{xy}$ such
that $\ol{w_i} = \lambda_i b(w_i)$. Hence, $B_{xy}=\{
b(w_1),\ldots, b(w_r)\}$ so  $\sigma = b(w_{i_0}) = \lambda_{i_0}
w_{i_0}$ for some $i_0$ such that $1\leq i_0 \leq r$. Then set
$p(\sigma)= w_{i_0}$. Moreover, if $p(\sigma)$ and $p'(\sigma)$
are two different such  paths and $\mu,\  \mu'$ scalars such that
$\sigma = \mu \ol{p(\sigma)} = \mu' \ol{p'(\sigma)}$ we get $\mu
\ol{p(\sigma)} - \mu' \ol{p'(\sigma)} = \sigma - \sigma = 0$ which
is a minimal relation, and this shows that $p(\sigma)$ and
$p'(\sigma)$ are naturally homotopic.

\item [$c)$] Let $w_1,w_2$ be two parallel homotopic paths, and suppose that $\ol{w}_1$ and $\ol{w}_2$
are linearly independent. Without loss of generality, we can
assume that  $w_1$ and $w_2$ appear  in the same minimal relation.
One then has  $\sum_{i=1}^n \lambda_i \ol{w}_i = 0$, where
$\lambda_i \in k^{*}$,  and the paths $w_i$ are parallel. Assume,
without loss of generality, that $\{\ol{w}_1, \ol{w}_2, \ldots,
\ol{w}_r \}$ is a maximal linearly independent set of $\{\ol{w}_1,
\ldots, \ol{w}_n \}$, and that,  after re-ordering if necessary
(recall that $A$ has a semi normed basis), there are scalars $a_j$
such that $\ol{w}_{r+1} = a_{r+1} \ol{w}_1,\ldots,$\
$\ol{w}_{r+i_1} = a_{r+i_1} \ol{w}_1,$ \ $\ol{w}_{r+i_1+1} =
a_{r+i_1+1} \ol{w}_2, \ldots $ $\ol{w}_n = a_n \ol{w}_r$. Thus,
replacing in $\sum_{i=1}^n \lambda_i \ol{w}_i = 0$ implies that
$$\lambda_1 + \lambda_{r+1} a_{r+1} + \cdots +
\lambda_{r+i_1}a_{r+i_1} = 0$$ and then $\lambda_1 \ol{w}_1
+\lambda_{r+1}\ol{w}_{r+1} + \cdots +
\lambda_{r+i_1}\ol{w}_{r+i_1} = 0$, which is a contradiction to
the minimality of the original relation.\cqfd
\end{enumerate}

\subsection{Theorem}\label{subsec:iso-complexes}{\em  Let $A=kQ/I$ be an algebra having a semi-normed basis, then:

\begin{enumerate}

    \item [$a)$] There exists an epimorphism of complexes
    \mor{{\rm SC}_\bullet (A)}{{\rm C}_\bullet ({\mathcal B}^\sharp)}{\phi_\bullet^\sharp},

    \item [$b)$] There exists an isomorphism of complexes
    \mor{{\rm SC}_\bullet (A)}{{\rm C}_\bullet ({\mathcal B})}{\phi_\bullet}.
\end{enumerate}}

\pf It is clear that ${\rm C}_0({\mathcal B}) \simeq {\rm C}_0({\mathcal B}^\sharp) \simeq {\rm SC}_0(A)$. For $n \geq 1$, consider the
morphisms $\phi^\sharp_n$ and $\phi_n$ defined by

$$\xymatrix@C=15pt@R=1pt{ \phi_n^\sharp:{\rm SC}_n(A) \ar[r]&{\rm
C}_n({\mathcal B}^\sharp )&\phi_n:{\rm SC}_n(A)\ar[r]&{\rm
C}_n({\mathcal B})\\ \hbox{ \ \ \  \ }(\sigma_1,\ldots, \sigma_n)
\ar@{|->}[r]& (\wt{p(\sigma_1)},\ldots, \wt{p(\sigma_n)}) &\hbox{
\ \ }(\sigma_1,\ldots, \sigma_n) \ar@{|->}[r]&
(\nwt{p(\sigma_1)},\ldots, \nwt{p(\sigma_n)}) }$$ It follows from
the second statement in lemma \ref{subsec:lemma_bases} that
$\phi_n^\sharp$ and $\phi_n$ are epimorphisms of abelian groups.
Moreover, it is clear that they commute with the boundary
operators involved. This proves the first statement. In order to
prove the second, consider the map $$\xymatrix@R=1pt{
 \psi_n:{\rm C}_n({\mathcal B}) \ar[r] & {\rm SC}_n(A)\\
\hbox{\ \ \ \ \ \ \ \ \ }(\nwt{\sigma}_1,\ldots, \nwt{\sigma}_n) \ar@{|->}[r]& (b(\sigma_1),\ldots, b(\sigma_n)).
}$$

The fact that the definition of $\psi_n$ is not ambiguous follows from the third statement in lemma \ref{subsec:lemma_bases}.  It is
clear that $\psi_n$  commutes with the boundary operators, so it defines a morphism of complexes. Finally, it is straightforward to check
that $\phi_\bullet$ and $\psi_\bullet$ are mutually inverse.\cqfd

\subsection{Corollary}\label{subsec:iso-homology}{\em With the above hypotheses,  for each $i\geq 0$, there are isomorphisms of
abelian groups
\begin{enumerate}

    \item [$a)$] ${\rm H}_i ({\mathcal B}) \simeq {\rm SH}_i (A)$, and

    \item [$b)$] ${\rm H}^i ({\mathcal B}, G) \simeq {\rm SH}^i (A, G)$, for every abelian group $G$.

\end{enumerate}}
\cqfd

\subsection*{Remark}In case $G$ is a commutative ring $R$, then the complexes involved are endowed with a canonical cup product. This
provides to ${\rm H}^\bullet (A,R) = \oplus_{i\geq 0} {\rm H}^i(A,
R)$, and ${\rm H}^\bullet ({\mathcal B}, R) = \oplus_{i\geq 0}
{\rm H}^i({\mathcal B}, R)$, a graded commutative ring structure.
A direct computation shows that the morphisms induced in
cohomology  $\phi_\bullet,\ \phi^\sharp_\bullet,$ and
$\psi_\bullet$ preserve these products. Thus, the morphisms in the
corollary above are ring homomorphisms.

\subsection*{Examples.}

\begin{enumerate}
\item Consider the quiver $Q$

$$\xymatrix@R=3pt{                   & 4\ar@/^/[dr]^{\beta_1}
&  &\\   6\ar@/^/[r]^{\alpha_1} \ar@/_/[r]_{\alpha_2} & 5
\ar@/_/[dr]_{\beta_3} \ar[r]^{\beta_2} & 2 \ar@/^/[r]^{\gamma_1}
\ar@/_/[r]_{\gamma_2}&1\\ &                   & 3 &}$$ \noindent
bound by $I=<(\alpha_1 - \alpha_2)\beta_3,\ \beta_1(\gamma_1 -
\gamma_2)>$. Even if in this case the complexes  ${\rm
SC}_\bullet(A)$ and ${\rm C}_\bullet({\mathcal B}^\sharp)$ are not
isomorphic, they still have the same homology groups. Let ${\rm
K}_\bullet = {\rm Ker}\ \phi_\bullet^\sharp$. We have that ${\rm
K}_n$ has a basis $\{(\sigma_1,\ldots,\sigma_n) -
(\sigma'_1,\ldots,\sigma'_n)|\ \sigma_i \sim \sigma'_i \hbox { for
all } i \hbox{ such that } 1\leq i\leq n \}$. More precisely,
${\rm rk\ K}_0 = 0,\ {\rm rk\ K}_1 = 7,\ {\rm rk\ K}_2 = 10,$ and
${\rm rk\ K}_3 = 3$. We leave to the reader the definition of a
contracting homotopy \mor{{\rm K}_\bullet}{{\rm K}_\bullet
[1]}{s_\bullet}.

\item Consider the quiver $Q$

$$\xymatrix{ x_1\ar@/^/[r]^{\alpha_1} \ar@/_/[r]_{\beta_1}&
x_2\ar@/^/[r]^{\alpha_2} \ar@/_/[r]_{\beta_2}&x_1}$$

\noindent bound by the ideal $I=<\alpha_1 \beta_2 + \beta_1
\beta_2 - \beta_1 \alpha_2,\ \alpha_1 \alpha_2 + \beta_1 \alpha_2
- \beta_1 \beta_2>$. The sets of cells of ${\mathcal B}$ are the
following: ${\mathcal C}_0 = \{x_1,x_2,x_3\}$, ${\mathcal C}_1 =
\{\nwt{\alpha_1},\nwt{\alpha_2}, \nwt{\beta_1}, \nwt{\beta_2},
\nwt{\alpha_1 \alpha_2}\}$, and ${\mathcal C}_2 =
\{(\nwt{\alpha_1}, \nwt{\alpha_2}),
(\nwt{\beta_1},\nwt{\alpha_2}), (\nwt{\alpha_1}, \nwt{\beta_2}),
(\nwt{\beta_1}, \nwt{\beta_2})\}$. On the other hand, the sets of
cells of ${\mathcal B}^\sharp$ are ${\mathcal C}^\sharp_0 =
\{x_1,x_2,x_3\}$, ${\mathcal C}^\sharp_1 = \{\wt{\alpha_1},
\wt{\alpha_2}, \wt{\alpha_1 \alpha_2}\}$, and ${\mathcal
C}^\sharp_2 = \{(\wt{\alpha_1}, \wt{\alpha_2})\}$.  An
straightforward computation yields to

$${\rm H}_i({\mathcal B}) =  \left\{ \begin{array}{rl}
 \mathbb{Z} & \mbox{ if } i \in\{0,2\} ,\\
 0      & \mbox{ otherwise. } \end{array} \right. \mbox{ and }
{\rm H}_i({\mathcal B}^\sharp) =  \left\{ \begin{array}{rl}
 \mathbb{Z} & \mbox{ if } i = 0 ,\\
 0      & \mbox{ otherwise. } \end{array} \right.
 $$

\noindent Thus, ${\mathcal B}$ and ${\mathcal B}^\sharp$ do not
have the
 same homotopy type. However, note that $A=kQ/I$, does not have a
 semi-normed basis.

\end{enumerate}

\medskip
Recall that, given an arcwise topological space $X$, the Hurewicz-Poincaré theorem (see \cite{ROT88}, for instance) states that
its first homology group, ${\rm H}_1(X)$ is the abelianisation of its fundamental group $\pi_1(X)$. This allows to prove
the following.

\subsection{Corollary}\label{subsec:constrainded}{\em Let $A=kQ/I$ be a constricted algebra having  a semi-normed basis. Then its
first simplicial (co)homology groups are independent of the presentation of $A$.}

\pf For such an algebra one has

$${\rm SH}_1(A)\simeq {\rm H}_1({\mathcal B}) \simeq \frac{\pi_1({\mathcal B})}{[\pi_1({\mathcal B}),\pi_1({\mathcal B})]} \simeq
\frac{\pi_1(Q,I)}{[\pi_1(Q,I),\pi_1(Q,I)]}.$$

On the other hand, the universal coefficients theorem gives

\begin{eqnarray*}
{\rm SH}^1(A,G) &\simeq &{\rm Hom}_{\mathbb Z}({\rm SH}_1(A), G) \oplus {\rm Ext}^1_{\mathbb Z}({\rm SH}_0(A), G)\\
        & \simeq &{\rm Hom}_{\mathbb Z}({\rm SH}_1(A), G).
        \end{eqnarray*}\cqfd

We now turn to the Hochschild cohomology of algebras.

%
%
%
%
%
%
%
%
%
%
%
%
%
%
%
%
%
%
%

\section{Hochschild cohomology}

Recall that for an arbitrary $k$-algebra $A$, its enveloping
algebra is the tensor product    $A^e = A\otimes_k A^{op}$. Thus,
an $A-A-$bimodule can be seen equivalently as an  $A^e$-module.
The Hochschild cohomology groups ${\rm HH}^i(A,M)$ of an algebra
$A$ with coefficients in some $A-A-$bimodule $M$, are, by
definition the groups ${\rm Ext}^i_{A^e}(A,M)$. In case $M$ is the
$A-A-$bimodule ${}_A A_A$, we simply denote them by ${\rm
HH}^i(A)$. We refer the reader to \cite{CE, Hap89, Redondo01}, for
instance, for general results about Hochschild (co)-homology of
algebras.

\subsection{A convenient resolution} In \cite{C3}, Cibils gave a convenient projective resolution of $A$ over $A^e$.
Let $E$ be the subalgebra of $A$ generated by the vertices of $Q$. Note that $E$ is semi-simple, and that
$A = \rad{A} \oplus E$ as $E-E-$bimodule. Let $\rad{A}^{\otimes n}$ denote the $n^{th}$ tensor power of $\rad{A}$ with
itself over $E$. With these notations, one has a projective resolution of $A$ as $A-A-$module:

$$\xymatrix@R=5pt{
\ldots \ar[r]& A \otimes_E \rad{A}^{\otimes n}\otimes_E A \ar[r]^{b_n}& A \otimes_E \rad{A}^{\otimes n}\otimes_E A \ar[r]^{\
\ \ \ \ \ \ \ \ \ \ \ \ \ b_{n-1}}&\ldots \\
       \ar[r]& A \otimes_E \rad{A} \otimes_E A \ar[r]^{b_1}             &A\otimes_E A \ar[r]^{b_0}       & A\ar[r] &0 }$$

where $d_0$ is the multiplication and
\begin{eqnarray*}
b_n(a\otimes \sigma_1 \otimes \cdots \otimes \sigma_n \otimes b) & = & a\sigma_1 \otimes \cdots \otimes \sigma_n \otimes b\\
&+& \sum_{j=1}^{n-1}(-1)^j a\otimes \sigma_1 \otimes \cdots \otimes \sigma_j\sigma_{j+1} \otimes \cdots \otimes \sigma_n b\\
&+& (-1)^n a\otimes \sigma_1 \otimes \cdots \otimes \sigma_n  b .
\end{eqnarray*}

Moreover, there is an obvious natural isomorphism
${\rm Hom}_{A^e}(A\otimes_E \rad{A}^{\otimes n } \otimes_E A, A) \simeq {\rm Hom}_{E^e}(\rad{A}^{\otimes n}, A)$. This will
be useful later. We denote by $b_n$ the corresponding boundary operator, and, moreover, we let
$b^n  = {\rm Hom}_{E^e}(b_n, A)$.

\subsection*{Remark} Note that the tensor products are taken over $E$. Thus, if $\sigma_1,\ \sigma_2 \in {\rm
rad}A$, with, say  $\sigma_1 \in e_i A e_j$, and  $\sigma_2 \in e_l A e_m$ then, in ${\rm rad}A^{\otimes 2}$, one has

$$\sigma_1 \otimes \sigma_2 = \sigma_1 e_j \otimes e_l \sigma_2 =  \sigma_1 \otimes  e_j  e_l \sigma_2$$ and this vanishes
if $j\not =l$. The same argument shows that  $\rad{A}^{\otimes n}$ is generated by elements of the form $\sigma_1 \otimes
\cdots \otimes \sigma_n$ where $\sigma_i \in e_{i-1} A e_i$ for $1 <i\leq n$. Moreover, if $A$ admits a semi-normed basis
$B$, then using $k-$linearity, one can assume that each $\sigma_i$ is an element of $B$.

\medskip
Following \cite{MdlP99}, for $n\geq 1$, define  \mor{{\rm Hom}_{\mathbb Z}({\rm SC}_n(A), k^+)}{{\rm
Hom}_{E^e}(\rad{A}^{\otimes n},A)}{\epsilon_n} in the following way:  for $f \in {\rm Hom}_{\mathbb Z}({\rm SC}_n(A),
k^+)$, and a basis element $(\sigma_1, \ldots \sigma_n)$, put   $\epsilon_n(f)(\sigma_1\otimes \cdots \otimes \sigma_n) =
f(\sigma_1,\ldots, \sigma_n) \sigma_1 \cdots \sigma_n$  whenever $(\sigma_1, \ldots, \sigma_n)$  $\in {\rm SC}_n(A)$,
and 0 otherwise.

Also, for $n\geq 1$, define  \mor{{\rm Hom}_{E^e}(\rad{A}^{\otimes n},A)}{{\rm Hom}_{\mathbb Z}({\rm SC}_n(A), k^+)}{\mu_n}
as follows: for a basis  element $(\sigma_1, \ldots, \sigma_n)$ in ${\rm SC}_n(A)$, we have  $\sigma_1 \sigma_2 \cdots
\sigma_n \not= 0$, and lies in, say, $e_0 A e_n$, which can be written as the direct sum of $k-$vector-spaces $<~\sigma_1
\sigma_2 \cdots \sigma_n > \oplus A'_{0n}$.  Moreover, $\sigma_1 \otimes \cdots \otimes \sigma_n \in {\rm rad}A^{\otimes
n}$, thus, for  $g \in{\rm Hom}_{E^e}(\rad{A}^{\otimes n},A)$ we have:

$$g(\sigma_1\otimes \cdots \otimes \sigma_n) = g(e_0 \sigma_1\otimes \cdots \otimes \sigma_n e_n) =e_0 g(\sigma_1\otimes
\cdots \otimes \sigma_n)e_n$$
\noindent so that $g(\sigma_1\otimes \cdots \otimes \sigma_n) \in e_0 A e_n$, and there is a scalar $\lambda$ and
$a_0\in A'_{0n}$ such that
$g(\sigma_1\otimes \cdots \otimes \sigma_n) = \lambda \sigma_1 \sigma_2 \cdots \sigma_n + a_0$. Define
$\mu_n(g)(\sigma_1, \ldots, \sigma_n ) = \lambda$.

\subsection{Lemma \cite{MdlP99}} \label{subsec:morphismes}{\em With the above notations, one has:

\begin{enumerate}

        \item [$a)$] $\mu_n \epsilon_n = id$, for $n\geq 1$,

        \item [$b)$] $\epsilon_{\bullet}$ is a morphism of complexes,

        \item [$c)$] If $A$ is schurian, then $\mu_{\bullet}$ is a morphism of complexes.

\end{enumerate} }\cqfd

\subsection*{Remark} In \cite{GS83}, it was shown that for incidence algebras,  the morphisms ${\rm H}^i(\epsilon)$ are isomorphisms.
Moreover, as a consequence of a result of \cite{PS01} (see also \cite{Bus02}),   ${\rm H}^1(\epsilon)$ is also an isomorphism for
schurian triangular algebras. Thus, in light of \cite{MdlP99} (or Lemma \ref{subsec:morphismes} above), one may naturally ask if in case
$A$ is schurian, the monomorphisms ${\rm H}^i(\epsilon)$ are isomorphisms. As the following example shows, this is not always the case.
\subsection*{Example} Consider the quiver

$$\xymatrix@R=4pt{          &2 \ar[dr]^{\beta}&  \\
1\ar[ur]^{\alpha} \ar[rr]_{\gamma}&     &3}$$

\noindent bound by the ideal $I = <\alpha \beta>$. The algebra $A=kQ/I$ is schurian, and one can easily compute
$${\rm SH}^i(A,k^+)=
\left \{
\begin{array}{ll}
k               & \mbox{if } i=0,1,\\
0               & \mbox{otherwise}
\end{array}
\right.
$$

\noindent On the other hand, using for instance, Happel's long exact sequence \cite{Hap89},  one gets

$$ {\rm HH}^i(A) =
\left \{
\begin{array}{ll}
k               & \mbox{if } i=0,1,2,\\
0               & \mbox{otherwise.}
\end{array}
\right.
$$

Recall that an algebra $A=kQ/I$ is said to be {\bf
semi-commutative} \cite{Hap89} whenever for $w, w'$ paths sharing
origin and terminus in $Q$, then $w \in I$ if and only if $w' \in
I$. For instance, incidence algebras are semi-commutative.
Schurian triangular, semi-commutative algebras are also called
{\bf weakly transitive} \cite{Dra}. The algebra in the preceding
example is schurian, but not semi-commutative. This leads to the
main result of this section, which is a generalization of a result
of Gerstenhaber and Schack \cite{GS83}, and makes more precise a
result of Martins and de la Pe\~na (theorem 3 in \cite{MdlP99}).
\subsection{Theorem} \label{subsec:cohomo-iso}{\em Let $A = kQ / I$ be a schurian triangular, semi-commutative algebra.
Then, for each $i \geq 0$, there is an isomorphism of abelian
groups $$\xymatrix{{\rm H}^i(\epsilon):{\rm SH}^i(A,k^+) \ar[r]^{\
\ \ \ \ \sim}&{\rm HH}^i(A)}.$$}

\pf In light of lemma \ref{subsec:morphismes}, there only remains to show that if $A$ is semi-commutative then $\epsilon_n
\mu_n = id$ for $n\geq 1$. Let $f \in {\rm Hom}_{E^e}(\rad{A}^{\otimes n},A)$, and $\sigma_1\otimes \cdots \otimes
\sigma_n$  be a basis element in ${\rm rad}A^{\otimes n}$, with, say $\sigma_1 \sigma_2 \cdots \sigma_n \in e_0 A e_n$.
Assume $\sigma_1 \sigma_2 \cdots \sigma_n \not= 0$. Since $A$ is schurian, there exists some scalar $\lambda$ such that
$f(\sigma_1\otimes \cdots \otimes \sigma_n) = \lambda \sigma_1\cdots \cdots \sigma_n$, thus  $(\mu_nf)(\sigma_1, \ldots,
\sigma_n) = \lambda$, and

\begin{eqnarray*} (\epsilon_n \mu_n f)(\sigma_1 \otimes \cdots \otimes \sigma_n)&=&((\mu_n f)(\sigma_1,
\ldots,\sigma_n))\sigma_1\cdots \sigma_n\\ &=& \lambda \sigma_1\ldots \sigma_n\\ &=& f(\sigma_1 \otimes \cdots \otimes
\sigma_n) \end{eqnarray*}

On the other hand, if $\sigma_1 \cdots \sigma_n = 0$, then,
$$(\epsilon_n \mu_n f)(\sigma_1 \otimes \cdots \otimes \sigma_n) = ((\mu_n f)(\sigma_1, \ldots,\sigma_n))\sigma_1\cdots \sigma_n\\
    = 0$$
Moreover, since $A$ is semi-commutative we have $e_0 A e_n = 0$, and therefore
$f(\sigma_1\otimes \cdots \otimes \sigma_n) = 0 $.
\cqfd

\subsection*{Remark.} Again, it is straightforward to check that $\epsilon_\bullet$ are $\mu_\bullet$ preserve cup-products, thus the
isomorphism above induce a ring isomorphism.

\subsection{Corollary \cite{GS83}}{\em  Let $(\Sigma,\leq)$ be a finite poset and $A=A(\Sigma)$ be its incidence algebra,
then, for each $i \geq 0$,  there is an isomorphism of abelian groups
${\rm H}^i({\mathcal  B}\Sigma, k^+) \simeq {\rm HH}^i(A)$.} \cqfd

\subsection*{Remark} There exist algebras which are not semi-commutative, but there are still an isomorphisms
$\xymatrix{\rm H^i(\epsilon): \ {\rm SH}^i(A,k^+) \ar[r] ^{\ \ \ \ \simeq}& {\rm HH}^i(A)}$ for all $i\geq 0$. Consider the following quiver $Q$
$$\xymatrix@R=2pt{              &4 \ar[r]^{\alpha_2}    &3 \ar[r]^{\alpha_3}& 2 \ar[dr]^{\alpha_4}& \\
 6 \ar[ur]^{\alpha_1} \ar[drr]_{\beta_1}&           &           &           &1  \\
                    &           &5 \ar[urr]_{\beta_2}&          &}$$and let
                    $A=kQ/I$ where $I=<\alpha_1 \alpha_2, \alpha_3\alpha_4>$. This algebra is schurian,
                    but not semi-commutative. One can easily compute
$${\rm SH}^i(A,k^+) = {\rm HH}^i(A) =
\left \{
\begin{array}{ll}
k               & \mbox{if } i=0,1,\\
0               & \mbox{otherwise.}
\end{array}
\right.
$$

Keeping in mind the last theorem, and proposition
\ref{subsec:retract}, we can get new short algebraic-topology
flavored proofs of some well-known results in \cite{Hap89,BM01}
about the Hochschild cohomology groups of monomial algebras. Let
$\chi(Q)$ be the Euler characteristic of $Q$, that is, let
$\chi(Q)=1-|Q_0| + |Q_1|$.

\subsection{Corollary \cite{Hap89}}{\em Let $A=kQ/I$ be a monomial
semi-commutative schurian algebra, then

\begin{enumerate}
    \item [$a)$] ${\rm HH}^0(A) = k$.

    \item [$b)$] ${\rm dim}_k {\rm HH}^1(A) = \chi(Q)$.

    \item [$c)$] ${\rm HH}^i(A) = 0$ for $i\geq 2.$
\end{enumerate} }

\pf With the above hypotheses, the graph $\ol{Q}$ is a strong deformation retract of ${\mathcal B}(Q,I)$, and theorem
\ref{subsec:cohomo-iso} holds. The results follow directly. \cqfd

\subsection{Corollary \cite{BM98}}{\em Let $A=kQ/I$ be a monomial algebra, then the following are equivalent:

\begin{enumerate}

    \item [$a)$] ${\rm HH}^i(A) = 0$ for $i>0.$

    \item [$b)$] ${\rm HH}^1(A) = 0.$

    \item [$c)$] $Q$ is a tree.

\end{enumerate}}

    \pf It is trivial that $a)$ implies $b)$. In order to show that $b)$ implies $c)$ assume that $Q$ is not a tree. Again, since $I$ is
    monomial, $\ol{Q}$ is a strong deformation retract of ${\mathcal B}(Q,I)$, and, since it is not a tree, we have ${\rm dim}_k {\rm
    SH}^1(A,k^+) = \chi(Q) > 0$. The result then follows from lemma \ref{subsec:morphismes}. Finally, we show that $c)$ implies $a)$. If
    $Q$ is a tree, then $A$ is schurian semi-commutative, thus theorem \ref{subsec:cohomo-iso} applies. But $\ol{Q}$, which is a strong
    deformation retract of  ${\mathcal B}(Q,I)$, is an acyclic $1-$dimensional complex. The result follows directly. \cqfd

\section*{Acknowledgements} We would like to thank professor
Ibrahim Assem for several helpful discussions and comments. Also,
S. Liu, who made an observation that led to distinguish ${\mathcal
B}$, and ${\mathcal B}^\sharp$ . The author gratefully acknowledge
financial support from the {\em Institut des Sciences
Mathématiques}, ISM.

\bibliography{biblio}
\end{document}